
\magnification=\magstep1
\input amstex
\documentstyle{amsppt}
\hsize = 6.5 truein
\vsize = 9 truein
\TagsAsMath
\NoRunningHeads
\topmatter
\title
A short proof of the irreflexivity conjecture
\endtitle
\author
Thomas Jech
\endauthor
\thanks Supported in part by NSF grant DMS--8918299
\endthanks
\address 
Department of Mathematics,
The Pennsylvania State University,
University Park, PA 16802, U.S.A.
\endaddress
\email jech\@math.psu.edu 
\endemail
\endtopmatter

\document
\baselineskip 18pt
 
Let $(A,*)$ be the free distributive algebra on one generator; $*$ is
a binary operation that satisfies
$$a*(b*c) = (a*b)*(a*c).$$
In his solution of the word problem for $(A,*)$ \cite{ }, Laver introduced 
the relation
$$a < b \qquad\text{ iff }\qquad
b=((a*b_1)*b_2*\cdots)*b_k \text{ for some } b_1,\dots,
b_k$$
and proved that $<$ is a linear ordering. Laver's proof uses a large cardinal 
assumption that is not provable in ordinary mathematics.
 
In a series of papers \cite{ }, Dehornoy showed (without any additional 
assumption) that any distinct $a,b\in A$ are 
comparable in $<$ and reduced the word 
problem to the problem of existence of a distributive one-generated
algebra $A$ in which $a < a$ does not hold for any $a \in A,$ the 
{\it irreflexivity conjecture.} A consequence of 
irreflexivity is that the algebra is free.
 
Dehornoy's latest paper \cite{ } proves the irreflexivity conjecture.
A major step of Dehornoy's proof is is the idea of embedding the free
distributive algebra into the braid group $B_\infty$ (see (5) below).
 
In this note we give a short proof of Dehornoy's theorem by giving
a simple argument showing that the operation $*$ in $B_\infty$
satisfies irreflexivity of $<.$
 
Let $F$ be the free group on free generators $x_1,x_2,x_3,\dots,$ and
let $B$ be the group of automorphisms of $F$ generated by
$\{\sigma_1,\sigma_2,\sigma_3,\dots\},$ where
$$
\sigma_i(x_i)=x_i x_{i+1} x_i^{-1},\quad\sigma(x_{i+1})=x_i,\quad
\sigma_i(x_j) = x_j \text{ if } i \neq j. \tag1
$$
 
\proclaim {Lemma 1} For every $i$ and every integer $k,$
$\sigma_i(x_i^k) = x_i x_{i+1}^k x_i^{-1}$ and $\sigma_i(x_{i+1}) =
x_i^k.$ \qed \endproclaim 
 
\proclaim {Lemma 2} (Artin \cite{ }.) For every $i$ and $j,$
$$
\sigma_i \sigma_{i+1} \sigma_i = \sigma_{i+1} \sigma_i \sigma{i+1},
\quad, \sigma_i \sigma_j = \sigma_j \sigma_i \text{ if }
|i-j| > 1. \tag2
$$\endproclaim
 
Let $W$ be the set of all (reduced) words in $F$ of the form
$$
w = x_{i_1}^{k_1}\cdots x_{i_n}^{k_n}, \tag3
$$
where $n \geq 1,$ $k_1,\dots,k_n \neq 0,$ and $i_1 \neq 1,$ $i_n \neq 1$
(words $\neq 1$ that do not begin or end with $x_1$ or $x_1^{-1}$).
 
\proclaim {Lemma 3} If $i \neq 1$ and $w \in W$ then $\sigma_i(w) \in W$
and $\sigma_i^{-1}(w) \in W.$ \endproclaim
 
\demo{Proof} Let $G$ be the subgroup of $F$ generated by $\{x_2,x_3,\dots\}$
and let $G^- = G - \{1\}.$ Clearly, $\sigma_i(x_1) = x_1,$ and if
$w \in G^-$ then $\sigma_i(w) \in G^-.$ If $w$ is in $W$, then either
$w \in G^-,$ or $w=u x_1^m v$ with $u,v \in G^-$ and $m \neq 0,$
or $w= u x_1^m z x_1^n v$ with $u,v \in G^-$ and $m,n \neq 0.$
In either case, $\sigma_i(w) \in W.$ The same holds for $\sigma_i^{-1}.$\qed
\enddemo
 
\proclaim{Lemma 4} If $w \in W,$ then $\sigma_1(x_1 w x_1^{-1}) =
x_1 \bar w x_1^{-1}$ for some $\bar w \in W.$\endproclaim
 
\demo{Proof} Let $F_2$ be the subgroup of $F$ generated by $\{x_1,x_2\}.$
 
{\bf Case I:} $w \in F_2.$ In this case,
$$
w = x_2^{k_1} x_1^{k_2} x_2^{k_3}\cdots x_1^{k_{n-1}} x_2^{k_n}
\qquad (n\geq 1,\, k_1,\dots,k_n \neq 0) \tag4
$$
and we have (by Lemma 1)
$$\align
\sigma_1(x_1 w x_1^{-1}) &=
x_1 x_2 x_1^{-1}\cdot x_1^{k_1}\cdot x_1 x_2^{k_2} x_1^{-1} \cdots
x_1^{k_n} \cdot x_1 x_2^{-1} x_1^{-1} \\
&=x_1 x_2 x_1^{k_1} x_2^{k_2}\cdots x_1^{k_n} x_2^{-1} x_1^{-1} \\
&= x_1 \bar w x_1^{-1}.
\endalign$$
{\bf Case II:} $w = u z v$ with $u,v \in F_2$ and $z \in Z$ where
$Z$ is the set of all words $ \neq 1$ that do not begin or end with
$x_1$, $x_1^{-1}$, $x_2$ or $x_2^{-1}.$ As $\sigma_1(z) \in Z,$
a similar argument as in Case I shows that $\sigma_1(x_1 w x_1^{-1})
= x_1 \bar w x_1^{-1}$ where $\bar w \in W:$ For instance consider the
alternative
$$
w = x_2^{k_1} x_1^{k_2} \cdots x_1^{k_{n-1}} x_2^{k_n} z
\qquad (n \geq 1,\, k_i \neq 0,\, z \in Z)
$$
(the other alternatives being similar). Then
$$
\sigma_1(x_1 w x_1^{-1}) = x_1 x_2 x_1^{k_1}\cdots x_2^{k_{n-1}}
x_1^{k_n-1}\cdot\bar z\cdot x_1 x_2^{-1} x_1^{-1} = x_1 \bar w x_1^{-1}.
\qed$$
\enddemo
 
\proclaim{Lemma 5} Let $\beta \in B$ be such that it can be written
as a word on $\{\sigma_1,\sigma_2,\dots\}$ with some occurrences
of $\sigma_1$ but no occurrences of $\sigma_1^{-1}.$ Then
$\beta(x_1) \neq x_1$ and so $\beta \neq 1.$
\endproclaim
 
\demo{Proof} This follows from Lemma 3, Lemma 4 and the fact that
$\sigma_1(x_1) = x_1 x_2 x_1^{-1} \in x_1 W x_1^{-1}.$ Since
$$
\beta = \alpha_1 \sigma_1 \alpha_2 \sigma_1 \cdots \sigma_1 \alpha_n,
$$
where each $\alpha_i$ is in the subgroup of $B$ generated by 
$\{\sigma_2,\sigma_3,\dots\},$ we have $\alpha_1(x_1) = x_1$,
$\alpha_1 \sigma_1(x_1) \in x_1 W x_1^{-1}$, and so
$\beta(x_1) \in x_1 W x_1^{-1}.$ \qed
\enddemo
 
{\bf Definition} (Dehornoy). For any $\alpha,\beta \in B,$ let
$$
\alpha * \beta = \alpha \cdot s(\beta) \cdot \sigma_1 \cdot s(\alpha^{-1})
\tag5$$
where $s$ is the shift: $s(\sigma_i) = \sigma_{i+1},$ $s(\sigma_i^{-1}) =
\sigma_{i+1}^{-1}$, and $s(\xi\cdot\eta) = s(\xi)\cdot s(\eta).$
 
\proclaim{Theorem} (Dehornoy \cite{ }.) For any $\alpha \in B$ and 
$\beta_1,\dots, \beta_k \in B,$
$$
\alpha \neq ((\alpha * \beta_1) * \beta_2 * \dots ) * \beta_k.
$$\endproclaim
 
\demo{Proof} It follows from (5) that $\alpha * \beta_1 * \dots *
\beta_k = \alpha\cdot\beta$ where $\beta$ has some occurrences of 
$\sigma_1$ but no occurrences of $\sigma_1^{-1}.$ Now apply Lemma 5.\qed
\enddemo
 
\proclaim{Lemma 6} (Dehornoy \cite{ }.) For all $\alpha,\beta,\gamma \in B,$
$\alpha *(\beta * \gamma) = (\alpha *\beta)*(\alpha*\gamma).$\endproclaim
 
\demo{Proof} Using (2).\qed\enddemo
 
\proclaim{Corollary} The free distributive algebra on one generator
satisfies the irreflexivity conjecture.\endproclaim

\Refs
\endRefs 
 
\enddocument
 
\bye